\documentclass{article}             
\usepackage{amssymb}

\def\bR{{\mathbb{R}}}

\def\fg{{\mathfrak{g}}}
\def\fh{{\mathfrak{h}}}
\def\fm{{\mathfrak{m}}}
\def\fn{{\mathfrak{n}}}

\newtheorem       {theorem}{Theorem}
\newtheorem{defin}[theorem]{Definition}
\newtheorem{prop} [theorem]{Proposition}
\newtheorem{lemma}[theorem]{Lemma}
\newtheorem{cor}  [theorem]{Corollary}

\begin{document}

\begin{center}
{\Large {The natural reductivity in Finsler geometry\\
in terms of geodesic graphs}}
\bigskip

{\large{Teresa Arias-Marco${}^1$ and Zden\v ek Du\v sek${}^2$} }
\bigskip

${}^1$ Department of Mathematics, University of Extremadura,\\
Av. de Elvas s/n, 06006 Badajoz, Spain\\
ariasmarco@unex.es\\
(ORCID: 0000-0003-0984-0367)

${}^2$ Institute of Technology and Business in \v Cesk\'e Bud\v ejovice\\
Okru\v zn\'\i\ 517/10, 370 01 \v Cesk\'e Bud\v ejovice, Czech Republic\\
zdusek@mail.vstecb.cz, corresponding author\\
(ORCID: 0000-0003-3073-9562)
\end {center}

\begin{abstract} 
A new geometrical definition of naturally reductive Finsler manifold
using geodeic graph is proposed, with a possible generalization.
Based on a construction from a recent paper by the authors,
Finsler metrics based on naturally reductive Riemannian metrics $g_i$ are studied.
Explicit examples of purely Finsler naturally reductive $\alpha_i$-type metrics are constructed.
Geodesic graphs on broad classes of Finsler $\alpha_i$-type metrics $F$
which are derived from naturally reductive Riemannian metrics and
which are not naturally reductive are described.
The influence of one-forms $\beta_j$ to the structure of geodesics of the metric $F$ is also demonstrated
and explicit construction of families of Finsler naturally reductive metrics
of the $(\alpha_i,\beta_j)$-type is described.
\end{abstract}
\bigskip
 
\noindent
{\bf MSClassification:} {53C22, 53C60, 53C30}
{\bf Funding:}
{The research is supported by grant GR24068
funded by Junta de Extremadura, partially funded by Fondo Europeo de Desarrollo Regional}

\section{Introduction}
If a connected Lie group $G$ acts on a Finsler manifold $(M,F)$ transitively by isometries,
the manifold is {\it homogeneous} and it can be 
naturally identified with the {\it homogeneous space} $(G/H,F)$, where $H$ is the isotropy subgroup.
If we denote by $\fg$ and $\fh$ the Lie algebras of $G$ and $H$, respectively,
we can always find a {\it reductive decomposition} $\fg=\fm+\fh$,
where $\fm\subset\fg$ is an ${\rm Ad}(H)$-invariant vector subspace
for the adjoint representation.
In a fixed reductive decomposition, the vector space $\fm$ is naturally identified with the tangent space $T_pM$,
for a fixed $p\in M$.
Using this identification, the Finsler metric $F$ on $M$ gives the ${\mathrm{Ad}}(H)$-invariant Minkowski norm
and the ${\mathrm{Ad}}(H)$-invariant fundamental tensor on $\fm$,
or, the ${\mathrm{Ad}}(H)$-invariant scalar product, respectively, in the case of a Riemannian metric.
A Riemannian homogeneous space $(G/H,g)$ is naturally reductive
if there exist a reductive decomposition $\fg=\fh+\fm$ such that
\begin{eqnarray}
\label{f1}
\langle [z,x]_\fm,y\rangle + \langle x,[z,y]_\fm\rangle = 0, \quad \forall x,y,z\in\fm.
\end{eqnarray}
It is well known that geodesics through the origin $p\in M=G/H$ of the Riemannian naturally reductive space
are orbits of particular one-parameter groups, namely they are of the form
\begin{eqnarray}
\label{f2}
\gamma(t) & = & {\mathrm{exp}}(ty)(p),\qquad y\in\fm,
\end{eqnarray}
with respect to the above decomposition.
For example, for symmetric spaces or normal homogeneous spaces,
the decomposition with the above property arise canonically.
In general, there is not a canonical reductive decomposition with the above property.
However, one can find it as the image of the geodesic graph.

Recently, there were several attempts to generalize the concept of natural reductivity
from Riemannian to Finsler geometry and following definitions were proposed.
The first definition generalizes the notion of a Minkowski Lie algebra.
\begin{defin}[\cite{La}]
\label{d1}
A Finsler homogeneous space $(G/H,F)$ is naturally reductive if there exist a reductive decomposition
$\fg=\fh+\fm$ such that
\begin{eqnarray}
\nonumber
g_y ( [x,u]_\fm,v ) + g_y ( u,[x,v]_\fm ) + 2 C_y ( [x,y]_\fm, u,v) = 0, \quad \forall x,y,u,v\in\fm, \quad y\neq 0.
\end{eqnarray}
\end{defin}
The second definition implicitly assumes that the Finsler metric $F$ is Berwald.
One of the possible definitions of a Berwald metric $F$ is that
its Chern connection is defined on $TM$ and it coincides with
the Levi-Civita connection of some Riemannian metric $g$.
\begin{defin}[\cite{DH}]
\label{d2}
A Finsler homogeneous space $(G/H,F)$ is naturally reductive if there exist a naturally reductive
Riemannian metric $g$ on $G/H$ such that the connections of $F$ and $g$ coincide.
\end{defin}
Here, we are missing the geometrical justification of the assumption that the given Finsler metric $F$
should be Berwald.
However, these definitions were shown to be equivalent in \cite{DYZ}.
We now propose a geometrical definition based on formula (\ref{f2}), which covers Definition \ref{d2}
and which is suitable for the study using geodesic graphs.
\begin{defin}
\label{d3}
A Finsler homogeneous space $(G/H,F)$ is naturally reductive if there exist a reductive decomposition
$\fg=\fh+\fm$ such that geodesics are of the form as in formula $(\ref{f2})$.
\end{defin}
If $(G/H,F)$ satisfies Definition \ref{d2}, then the Riemannian metric $g$ admits a reductive decomposition
with geodesics of the form (\ref{f2}).
Because connections of $F$ and $g$ coincide, their geodesics are identical and the same reductive decomposition
satisfies Definition \ref{d3}.

Natural reductivity is a special case of a more general property.
A homogeneous Finsler space $(G/H, F)$
is called a {\it geodesic orbit space} or just {\it g.o. space},
if every geodesic $\gamma(t)$ through $p\in G/H$ is
an orbit of a one-parameter subgroup of the group of isometries $G$,
namely it is of the form
\begin{eqnarray}
\label{f3}
\gamma(t) & = & {\rm exp}(ty)(p), \qquad y\in\fg.
\end{eqnarray}
Vectors which generate geodesics as in formulas (\ref{f2}) or (\ref{f3}) are called {\it geodesic vectors}.
We should also remark that the g.o. property, as well as the natural reductivity of the manifold $M$,
depends on the choice of the transitive isometry group $G$.
Naturally reductive manifolds and g.o. manifolds were studied broadly in Riemannian geometry by many authors,
for results and references, see for example the monograph \cite{BN} by V. Berestovskii and Yu. Nikonorov
and the references therein.

One of the techniques for the study of g.o. spaces is based on geodesic graphs
proposed originally by J. Szenthe in \cite{Sz} in the affine setting.
This concept was further developed by O. Kowalski and S. Nik\v cevi\'c in \cite{KNi} for the Riemannian setting.
\begin{defin}
\label{d4}
Let $(G/H,F)$ be a Finsler g.o. space and $\fg=\fh+\fm$ be a reductive decomposition.
A geodesic graph is an ${\mathrm{Ad}}(H)$-equivariant map $\xi\colon\fm\rightarrow\fh$
such that $y +\xi(y)$ is a geodesic vector for each $o\neq y\in\fm$.
\end{defin}
In a g.o. space, at least one geodesic graph must exist.
In general, the components of a Riemannian geodesic graph are rational functions $\xi_i=P_i/P$,
where $P_i$ and $P$ are homogeneous polynomials and ${\mathrm{deg}}(P_i)={\mathrm{deg}}(P)+1$.
The existence of a linear geodesic graph is equivalent with the natural reductivity of the space $G/H$
and the image of the linear geodesic graph is the linear subspace which satisfies formula (\ref{f1}).
In the reductive decomposition which satisfies formula (\ref{f1}), geodesic graph is the zero map.
See for example the survey paper \cite{DuS2} by the second author for
details and references to geodesic graphs on various classes of Riemannian g.o. manifolds.
Definition \ref{d3} together with Definition \ref{d4} are in particular usefull for disproving
that a given metric $F$ is naturally reductive, in the situation when the unique geodesic graph is nonlinear.

In Finsler geometry, g.o. spaces started to be studied recently.
Relations of some Finsler g.o. spaces with Riemannian g.o. spaces
for Randers g.o. metrics on spheres $S^{2n+1}={\mathrm{U}}(n+1)/{\mathrm{U}}(n)$
and for some weakly symmetric metrics on nilpotent groups
were investigated by Z. Yan and S. Deng in \cite{YD}.
Geodesic orbit Finsler metrics on spheres were studied also by M. Xu in \cite{Xu}
and geodesic orbit $(\alpha,\beta)$ metrics were studied in \cite{Y} by Z. Yan.
Naturally reductive Finsler spaces were studied in the papers
\cite{DYZ}, \cite{DH} which we have already mentioned and some others.
In \cite{TX}, it was shown that $(\alpha_1,\alpha_2)$ metrics are naturally reductive
if and only if they are an $f$-product of two naturally reductive Riemannian metrics.
Unfortunately, none of the papers about the topic contained examples of purely Finsler
naturally reductive spaces.

The investigation of geodesic graphs on Finsler g.o. spaces started with the papers
\cite{DuCMUC} and \cite{DuWS} by the second author, for Randers geodesic orbit metrics
and special families of weakly symmetric Finsler metrics on modified H-type groups.
New structures of geodesic graphs were observed and also relations
of Finslerian geodesic graphs with the Riemannian geodesic graphs.
In \cite{Du}, second author studied general homogeneous Finsler $(\alpha,\beta)$ metrics.
It was proved that all $(\alpha,\beta)$ metrics arising from a Riemannian g.o. metric $\alpha$
and an invariant one-form $\beta$ are Finsler g.o. metrics, possibly with respect to an extended isometry group.
Relations of Finslerian geodesic graphs with the Riemannian geodesic graphs for these metrics were described.
Examples of these situations were shown in
\cite{AMD}, where the present authors constructed geodesic graphs
for geodesic orbit $(\alpha,\beta)$ metrics on spheres and classified
these metrics on spheres and on projective spaces.

In the paper \cite{AMD2}, the present authors investigated g.o. property of special metrics
from the broad class of Finsler metrics which were proposed in \cite{JS} and
whose Minkowski norms $F=F(y)$ on $\fm$ arise as
\begin{eqnarray}
\label{ma}
F & = & \sqrt{L(F_1,\dots, F_k, \beta_1,\dots, \beta_l) },
\end{eqnarray}
where $F_i=F_i(y)$ are Minkowski norms and $\beta_j=\beta_j(y)$ are one-forms on $\fm$.
We analyzed the g.o. property of Finsler metrics of the type (\ref{ma}), without one-forms
and with the assumption that $F_i$ are norms determined by Riemannian g.o. metrics $g_i$, hence $F_i=\sqrt{g_i}$.
With the natural assumption that the Riemannian metrics $g_i$ are positively related,
we have proved that the invariant Finsler metric corresponding to the Minkowski norm
\begin{eqnarray}
\label{ma2}
F & = & \sqrt{L(\sqrt{g_1},\dots, \sqrt{g_k} ) }
\end{eqnarray}
is also g.o. metric and we described the Finslerian geodesic graph.
An important observation shows that the assumption of positively related metrics $g_i$ above
implies that the metrics (\ref{ma2}) are of the $\alpha_i$-type, which includes $(\alpha_1,\alpha_2)$ type in particular.
If it is constructed on the product $M=M_1\times M_2$, it includes the $f$-product metrics.

In the present paper, in Section \ref{S3}, we focus on geodesic graphs for metrics of the type (\ref{ma2}),
where $F_i$ are Minkowski norms of naturally reductive Riemannian metrics.
We describe geodesic graphs of these metrics,
under the assumption that the initial naturally reductive Riemannian metrics are positively related.
If the construction is done on the product manifold $M=M_1\times M_2$,
we obtain the first explicit examples of naturally reductive purely Finsler $(\alpha_1,\alpha_2)$ metrics,
or, in general, $\alpha_i$-type metrics.
If the manifold $M$ is irreducible, we obtain the nonlinear geodesic graph and the manifold
with the Finsler metric is g.o., but not naturally reductive.
We illustrate the construction with an example of the Heisenberg group $H_3$.
If the assumption that the initial naturally reductive metrics $g_i$ are positively related is not valid,
one has to use the more general approach using nonlinear geodesic graphs,
described in \cite{AMD2}.
As an example of a family of naturally reductive Riemannian metrics which are not positively related,
one can consider the family of metrics on the sphere $S^7$, described in detail also in \cite{AMD2}.

In Section \ref{S4}, we also study the special case of Finsler metrics of the type $F=\sqrt{L(\sqrt{g},\beta)}$,
where $g$ is a naturally reductive Riemannian metric.
In Finsler geometry, the Berwald type assumption in Definition \ref{d2} on the metric $F$ is very strong and limiting.
For example, it is well known that an $(\alpha,\beta)$ metric $F$ is Berwald if and only if the one-form $\beta$
is parallel with respect to the Levi-Civita connection $\nabla$ of $\alpha$.
But the existence of a parallel one-form, or, equivalently, the existence of a parallel vector field,
implies that the manifold is locally reducible.
From the available literature, we observed general expectations that $(\alpha,\beta)$ Berwald metrics $F$
are naturally reductive if and only if the Riemannian metric $\alpha$ is naturally reductive.
In such a case, the Chern connection of $F$ is the same as the Levi-Civita connection of $\alpha$.
However, we are not aware of a clear and simple proof of this fact.
Concerning Finsler metrics of the type $F=\sqrt{L(\sqrt{g},\beta)}$,
we demonstrate easily, using geodesic graphs, that the Levi-Civita connection of $g$ and
the Chern connection of $F$ share the same geodesics and consequently they must coincide
if and only if $\beta$ is parallel with respect to the Levi-Civita connection of $g$.
Finally, we will provide explicit construction
of examples of naturally reductive $(\alpha_i,\beta)$-type Finsler spaces.

\section{Preliminaries}
Recall that we consider a homogeneous space $M=G/H$ with a fixed reductive decomposition $\fg=\fh+\fm$.
Let $p\in M$ and $F_p$ be an ${\mathrm{Ad}}(H)$-invariant Minkowski norm on $\fm\simeq T_pM$.
The Minkowski norm $F_q$ at any $T_qM, q\in M$ is obtained by the formula
\begin{eqnarray}
\nonumber
F_q(\sigma_*X) = F_p(X), \qquad X\in T_pM,
\end{eqnarray}
where $\sigma\in G$ and $\sigma(p)=q$.
It gives the correspondence of invariant Finsler metrics $F$ on $M$ and the ${\mathrm{Ad}}(H)$-invariant
Minkowski norms $F_p$ on $\fm$.
We shall work just with the vector space $\fm$, the ${\mathrm{Ad}}(H)$-invariant Minkowski norm
on $\fm$ and its fundamental tensor, which we denote also by $F$ and $g$, respectively.
We recall first the definition of a homogeneous $\alpha_i$-type metric, which covers
homogeneous $(\alpha_1,\alpha_2)$ metrics and also the $f$-product metrics if $M=M_1\times M_2$.
\begin{defin}[\cite{AMD2}]
\label{ai}
Let $(G/H,F)$ be a homogeneous Finsler space with a reductive decomposition $\fg=\fm+\fh$.
Consider the irreducible decomposition
\begin{eqnarray}
\nonumber
\fm & = & \oplus_{i=1}^s \fm_i.
\end{eqnarray}
with respect to the adjoint action of $H$ on $\fm$.
Choose symmetric positively definite ${\mathrm{Ad}}(H)$-invariant bilinear forms $\alpha_i$ on $\fm_i, i=1\dots s$
and denote by $y_i$ the corresponding projections of a vector $y\in\fm$ onto $\fm_i$.
The Minkowski norm $F$ on $\fm$ and the corresponding homogeneous Finsler metric on $G/H$ is of the $\alpha_i$-type
if there exist a smooth function $f\colon [ 0,\infty)^s\rightarrow{\mathbb{R}}$ such that
\begin{eqnarray}
\label{f5}
F(y) & = & f(\alpha_1(y_1),\dots, \alpha_s(y_s)),\qquad y\in\fm.
\end{eqnarray}
\end{defin}
An unpleasant feature of this definition is that it does not give sufficient conditions for the function $f$.
Clearly, the smoothness is not a sufficient condition, because the functions $\alpha_i$ are not
injective and also the definition domain of the function $f$ is not an open set.
Hence, at points where some of the $\alpha_i$'s reach the zero value, even the smoothness of the function $F$
is not guaranteed. According to our knowledge, the explicit examples of such metrics can be constructed
using the notion of positively related metrics introduced in \cite{AMD2} and the approach proposed in \cite{JS}.

\begin{defin}[\cite{AMD2}]
\label{galpha}
Let $G/H$ be a homogeneous space with a reductive decomposition $\fg=\fm+\fh$,
the ${\mathrm{Ad}}(H)$-invariant irreducible decomposition $\fm=\oplus_{i=1}^s \fm_i$
and initial ${\mathrm{Ad}}(H)$-invariant scalar products $\alpha_i$ on the respective spaces $\fm_i$.
Consider the family of scalar products
\begin{eqnarray}
\label{alpha}
g_{c^1 \dots c^s} & = & \sum_{i=1}^s c^i\cdot \alpha_i,
\end{eqnarray}
for any numbers $0<c^i\in {\mathbb{R}}$, and related family of invariant Riemannian metrics on $G/H$.
This family of scalar products on $\fm$ and corresponding family of invariant Riemannian metrics on $G/H$
will be called scalar products positively related and metrics positively related.
\end{defin}
In \cite{JS}, the authors consider conic pseudo-Finsler metrics $F$ which arise
from given Finsler metrics and one-forms and they are given by the formula
\begin{eqnarray}
\label{ma3}
F^2(y) & = & L(F_1,\dots, F_k, \beta_1,\dots, \beta_l),
\end{eqnarray}
for a continuous function $L$, which is smooth and positive away from $0$ and positively homogeneous of degree 2.
According to Theorem 4.1 in \cite{JS}, to obtain a Minkowski norm $F$, the function $L$ must satisfy:\\
(i) $L_{,i} \geq 0$, for $i=1,\dots,k$,\\
(ii) ${\mathrm{Hess}}(L)$ be positive semi-definite,\\
(iii) $L_{,1} + \dots + L_{,k} >  0$.\\
Here the comma in $L_{,i}$ means the derivative with respect to the corresponding coordinate.
We will be interested only in positively definite Finsler metrics defined on whole $\fm$.
However, to keep simplicity, we will not determine explicitly conditions for the definition domain to be all of $\fm$.
If there is just one Riemannian metric $F_1$ and one one-form $\beta_1$, 
metrics of the type $F(y) = \sqrt{L( \sqrt{g}, \beta )}$ give an alternative approach
to $(\alpha,\beta)$ metrics, see \cite{JS}.

If there are no one-forms $\beta_j$ in formula (\ref{ma3}) and the Finsler metrics $F_i$ are Riemannian
and positively related, then the Finsler metric $F$ can be written also using formula (\ref{f5}),
for some function $f$ determined by $L$, and hence the Finsler metric $F$ is of the $\alpha_i$-type.
For an illustrating example, we will consider later the metrics of the type
\begin{eqnarray}
\label{Fex}
F(y) & = & \left ( \sum_{k=1}^r F_k(y)^q \right )^{\frac{1}{q}}, \qquad q>2,
\end{eqnarray}
for Riemannian naturally reductive metrics $F_k=\sqrt{g_k}$.

Throughout the rest of the paper, we always assume that the smooth and positively homogeneous function $L$
satisfies the conditions (i)-(iii) above.
Let us now recall the geodesic lemma, for a general Finsler metric $F$.
\begin{lemma}[\cite{La}]
\label{golema2}
Let $(G/H,F)$ be a homogeneous Finsler space with a reductive decomposition $\fg=\fm+\fh$.
A nonzero vector $y\in{\fg}$ is geodesic vector if and only if it holds
\begin{equation}
\label{gl2}
g_{y_\fm} ( y_{\mathfrak m}, [y,u]_{\mathfrak m} ) = 0 \qquad \forall u\in{\mathfrak m},
\end{equation}
where the subscript $\fm$ indicates the projection of a vector from $\fg$ to $\fm$.
\end{lemma}
To adapt geodesic lemma for our type of Finsler metrics,
we will need the technical lemma to express the fundamental tensor.
The special case of this lemma was proved in \cite{AMD2}, without the one-forms. See also \cite{JS}.

\begin{lemma}
Let $g_1,\dots, g_k$ be homogeneous Riemannian metrics, $\beta_1,\dots,\beta_l$ be homogeneous one-forms on $G/H$
and let $\fg=\fm+\fh$ be a reductive decomposition.
We use the same notation for the bilinear forms and one-forms on $\fm$.
Let $F=\sqrt{L(\sqrt{g_1},\dots,\sqrt{g_k},\beta_1,\dots,\beta_l)}$ on $\fm$, which gives a homogeneous Finsler metric on $G/H$.
For arbitrary vectors $y,v\in\fm$, the fundamental tensor $g$ of $F$ satisfies the formula
\begin{eqnarray}
\label{fg}
2 g_y(y,v) & = & \sum_{j=1}^k \frac{ L_{,j}}{F_j(y)} \cdot g_j(y,v) + \sum_{m=k+1}^{k+l} L_{,m}\cdot \beta_m(v).
\end{eqnarray}
\end{lemma}
{\it{Proof.}}
Let $y,v\in\fm$ be arbitrary fixed vectors.
For the vector $w=y+tv\in\fm$, where $t\in\bR$, we have
$F^2(w) = L(\sqrt{g_1(w,w)},\dots, \sqrt{g_k(w,w)},\beta_1(w),\dots,\beta_l(w))$.
The result is obtained by the differentiation of $F^2(w)$ with respect to $t$, using the chain rule and putting $t=0$.
$\hfill\square$

We now assume that the metrics $g_j$ are positively related.
Decompose the scalar products $g_j$ according to formula (\ref{alpha}) as
$g_j = \sum_{i=1}^s c^i_j\cdot \alpha_i$ and denote
\begin{eqnarray}
\label{F12}
B_j(y) =  \frac{ L_{,j}}{F_j(y)}, \qquad
C_i(y) =  \sum_{j=1}^k B_j(y) \cdot c^i_j = \sum_{j=1}^k{ \frac{ L_{,j}}{F_j(y)} } \cdot c^i_j.
\end{eqnarray}
Geodesic lemma for this situation is in the following form.
\begin{lemma}
\label{golemanew}
Let $G/H$ be a homogeneous space with a reductive decomposition $\fg=\fm+\fh$
and with the ${\mathrm{Ad}}(H)$-irreducible decomposition $\fm=\oplus_{i=1}^s \fm_i$.
Let the metrics $g_j$ be positively related with corresponding decompositions $g_j=\sum_{i=1}^s c_j^i \alpha_i$.
Let $F=\sqrt{L(\sqrt{g_1},\dots,\sqrt{g_k},\beta_1,\dots,\beta_l)}$ on $\fm$, which gives a homogeneous Finsler metric on $G/H$.
The vector $y+\xi(y)$, where $y\in\fm$ and $\xi(y)\in\fh$,
is geodesic vector for the Finsler metric $F$ if and only if for any $u\in\fm$ it holds
\begin{eqnarray}
\label{golemaf}
\sum_{i=1}^s C_i(y) \cdot \alpha_i\Bigl (y,[y+\xi(y),u]_{\mathfrak m}\Bigr )
+ \sum_{m=k+1}^{k+l} L_{,m}\cdot \beta_m( [y+\xi(y),u]_{\mathfrak m} ) & = & 0.\qquad
\end{eqnarray}
\end{lemma}
{\it{Proof.}}
In geodesic lemma, use formula (\ref{fg}) with $v=[y+\xi(y),u]_\fm$,
decompose $g_j$ into $\alpha_i$ and use the definition of the functions $C_i(y)$.
$\hfill\square$

\section{Example: Finsler $\alpha_i$-type metrics on $H_3$}
\label{sech}
We are going to illustrate the concepts from Preliminaries on the $3$-dimensional Heisenberg group,
whose ${\mathrm{Ad}}(H)$-invariant decomposition is $\fm=\fm_1\oplus\fm_2$ and it admits a 2-parameter
family of positively related naturally reductive Riemannian metrics.
We will consider Finsler metrics of the form (\ref{Fex}) with $q=3$, which lead to $(\alpha_1,\alpha_2)$ Finsler metrics.
We apply geodesic lemma and we describe in detail the geodesic graph, which will be the topic of the further general
discussion.

The Lie algebra $\fn$ of the $3$-dimensional Heisenberg group $H_3$ is generated by the basis $B=\{E_1,E_2,E_3\}$
with the nontrivial Lie bracket
\begin{eqnarray}
\nonumber
[E_1,E_2] = E_3.
\end{eqnarray}
The $2$-parameter family of naturally reductive Riemannian metrics on $H_3$ is determined, up to a scalar multiple,
by the scalar product in which the above basis $B$ is orthogonal and
\begin{eqnarray}
\label{scalar}
\langle E_1,E_1 \rangle = \langle E_2,E_2 \rangle = 1, \quad \langle E_3,E_3 \rangle = c
\end{eqnarray}
for $c>0$.
We shall verify this fact by the construction of the linear geodesic graph. We express $H_3$ as a homogeneous space
$G/H$, where $H={\mathrm{SO}}(2)$ and $G=H_3\rtimes{\mathrm{SO}}(2)$.
The Lie algebra $\fh$ is generated by the operator $D$ with the nontrivial action
\begin{eqnarray}
\nonumber
D(E_1)=E_2, \quad D(E_2)=-E_1.
\end{eqnarray}
The subspaces $\fm_i$ in the ${\mathrm{Ad}}(H)$-invariant decomposition $\fn=\fm=\fm_1\oplus\fm_2$ above are
obviously $\fm_1={\mathrm{span}}(E_1,E_2)$ and $\fm_2={\mathrm{span}}(E_3)$.
The scalar product above can be obviously written as $\alpha_1 + c\cdot \alpha_2$, for the
basic scalar products $\alpha_i$ on $\fm_i$.
We substitute now the vectors $y = y_1 E_1 + y_2 E_2 + y_3 E_3$ and $\xi(y) = \xi_1 D$
into the Riemannian geodesic lemma
\begin{eqnarray}
\nonumber
g( y ,[ y + \xi(y) ,u]) & = & 0
\end{eqnarray}
and for $u$ we substitute, step by step, the elements of the above basis $B$.
We obtain the system of equations
\begin{eqnarray}
\xi_1 y_2 & = & cy_2 y_3, \cr
-\xi_1 y_1 & = & -cy_1 y_3,
\end{eqnarray}
for the component $\xi_1$, depending on components $y_j$.
Obviously, there is the unique solution
\begin{eqnarray}
\xi_1 = c y_3.
\end{eqnarray}
Let us consider now two naturally reductive metrics $g_1$ and $g_2$ of the given type, where we put
$c=c_1$ for $g_1$ and $c=c_2\neq c_1$ for $g_2$ in formula (\ref{scalar}).
Obviously, the metrics $g_i$ are positively related and any Finsler metric $F=\sqrt{L(\sqrt{g_1},\sqrt{g_2})}$
for admissible function $L$ described in Preliminaries gives a Finsler $(\alpha_1,\alpha_2)$ metric.
For the illustration, we consider the Finsler metric
\begin{eqnarray}
\label{metric}
F(y) & = & \Bigl ( \sum_{i=1}^2 F_i(y)^q\Bigr )^{\frac{1}{q}},
\end{eqnarray}
in the particular case $q=3$ with two Finsler metrics $F_i(y)=\sqrt{g_i(y,y)}$, for Riemannian metrics $g_i$ above.
The straightforward calculations give us
\begin{eqnarray}
L & = & F^2 = \Bigl ( \sum_{i=1}^2 F_i(y)^3\Bigr )^{\frac{2}{3}}, \cr
L_{,k} & = & \frac{2}{3} \Bigl ( \sum_{i=1}^2 F_i(y)^3\Bigr )^{-\frac{1}{3}}\cdot 3 \cdot {F_k}^2 = \frac{{2F_k}^2}{F}, \cr
B_{k} & = & \frac{L_{,k}}{2 F_k} = \frac{F_k}{F}.
\end{eqnarray}
Now, let us express the two metrics $g_i$ above as
$g_1=\alpha_1 + c_1\alpha_2$ and $g_2=\alpha_1 + c_2\alpha_2$.
We easily obtain, using notation of formula (\ref{F12}),
\begin{eqnarray}
\nonumber
2\, g_y(y,v) & = & B_1(y) g_1(y,v) + B_2(y) g_2(y,v) = \cr
         & = & B_1(y) \Bigr ( \alpha_1(y,v) + c_1 \alpha_2(y,v) \Bigr ) 
             + B_2(y) \Bigr ( \alpha_1(y,v) + c_2 \alpha_2(y,v) \Bigr ) = \cr
         & = & \Bigl ( B_1(y) + B_2(y) \Bigr ) \alpha_1(y,v) + \Bigl ( c_1 B_1(y) + c_2 B_2(y) \Bigr ) \alpha_2(y,v) = \cr
         & = & \Bigl ( \frac{ F_1 + F_2}{F} \Bigr ) \alpha_1(y,v) + \Bigl ( \frac{c_1 F_1 + c_2 F_2}{F} \Bigr ) \alpha_2(y,v) = \cr
         & = & C_1(y) \alpha_1(y,v) + C_2(y) \alpha_2(y,v).
\end{eqnarray}
Now, geodesic lemma with our metric $F(y)$ in formula (\ref{metric}) for $q=3$ and with the same setting
as above gives us the system of equations
\begin{eqnarray}
\xi_1 C_1(y) y_2 & = & C_2(y) y_2 y_3, \cr
-\xi_1 C_1(y) y_1 & = & -C_2(y) y_1 y_3,
\end{eqnarray}
with the unique solution
\begin{eqnarray}
\nonumber
\xi_1 = \frac{C_2(y)}{C_1(y)} y_3 =
\frac{c_1 F_1 + c_2 F_2} { F_1 +  F_2} y_3 =
\frac{c_1 \sqrt{ y_1^2 + y_2^2 + c_1y_3^2 } + c_2 \sqrt{ y_1^2 + y_2^2 + c_2y_3^2 } }
     { \sqrt{ y_1^2 + y_2^2 + c_1y_3^2 } +  \sqrt{ y_1^2 + y_2^2 + c_2y_3^2 } }  y_3.
\end{eqnarray}

\begin{cor}
\label{c2}
There exist a family of Finsler $(\alpha_1,\alpha_2)$ metrics $F$ on $H_3$
which arise as $F=\sqrt{L(\sqrt{g_1},\sqrt{g_2})}$, where $g_i$ are naturally reductive Riemannian metrics.
These Finsler metrics are g.o., but not naturally reductive.
\end{cor}

\section{Geodesic graph for metrics $F=\sqrt{L(\sqrt{g_1},\dots,\sqrt{g_k})}$}
\label{S3}
We are going to describe the general structure of geodesic graphs for Finsler metrics of the type
$F=\sqrt{L(\sqrt{g_1},\dots,\sqrt{g_k})}$, where $g_1,\dots,g_k$ are metrics
from the family of positively related naturally reductive Riemannian metrics.
The behaviour was illustrated above with the example of $H_3$.

Let $M=G/H$ be an irreducible homogeneous space with the ${\mathrm{Ad}}(H)$-irreducible decomposition
$\fm=\oplus_{i=1}^s\fm_i$ with at least two summands.
Let $g_{c^1\dots c^s}, c^i>0$, be a family of positively related naturally reductive Riemannian metrics on $M$,
with decompositions $g_{c^1\dots c^s}=\sum_{i=1}^s c^i\alpha_i$, corresponding to the decomposition of $\fm$.
Let $g_1,\dots,g_k$ be metrics from the above family.
We have observed in Preliminaries, that any Finsler metric $F=\sqrt{L(\sqrt{g_1},\dots,\sqrt{g_k})}$
for an admissible function $L$, is of the $\alpha_i$-type.

Let us further observe the general structure of linear geodesic graph for naturally reductive Riemannian metrics
$g_{c^1\dots c^s}$.
Geodesic lema for any of these metrics is in the form
\begin{eqnarray}
\sum_{i=1}^s c^i \cdot \alpha_i\Bigl (y,[y+\xi(y),u]_{\mathfrak m}\Bigr ) & = & 0, \quad \forall u\in\fm.
\end{eqnarray}
For given $u\in\fm_k$, it holds $[\xi(y),u]\in\fm_k$ and we can rewrite it in the form
\begin{eqnarray}
\alpha_k\Bigl (y,[\xi(y),u]_{\mathfrak m}\Bigr ) & = &
- \sum_{i=1}^s \frac{c^i}{c^k} \cdot \alpha_i\Bigl (y,[y,u]_{\mathfrak m}\Bigr ).
\end{eqnarray}
For each $u\in\fm$, on the left-hand side, we have linear terms, where variables are
coordinates $y_r$ of the vector $y$
and coefficients are scalar multiples of components $\xi^j$ of the vector $\xi(y)$.
On the right-hand side, we have quadratic terms in coordinates $y_r$
whose coefficients are rational functions in $c^i$.
Hence, any linear solution $\xi(y)$ of this system has the form as in formula (\ref{ft2}) below.
\begin{theorem}
\label{t1}
Let the components of the linear geodesic graph of the individual naturally reductive
Riemannian metrics $g_{c^1\dots c^s}$ be
\begin{eqnarray}
\label{ft2}
\xi_{c^1\dots c^s}^j & = & \sum_{r=1}^n k^j_r(c^1,\dots,c^s) y_r, \qquad j=1,\dots,{\mathrm{dim}}(\fh),
\end{eqnarray}
where $k^j_r(c^1,\dots,c^s)$ are fixed rational functions of the parameters $c^i$ determining the individual
metrics $g_{c^1\dots c^s}$.
Then the Finsler metric $F=\sqrt{L(\sqrt{g_1},\dots,\sqrt{g_k})}$ is a g.o. metric whose geodesic graph components
are given by the formula
\begin{eqnarray}
\label{ft3}
\xi^j & = & \sum_{r=1}^n k^j_r \left ( C_1(y),\dots,C_s(y)\right ) y_r, \qquad j=1,\dots,{\mathrm{dim}}(\fh),
\end{eqnarray}
where the functions $C_i(y)$ are given by formula $(\ref{F12})$.
\end{theorem}
{\it{Proof.}}
According to Lemma \ref{golemanew}, geodesic lemma for the Finsler metric $F$ for given $u\in\fm_k$ becomes
\begin{eqnarray}
\label{f22}
\alpha_k\Bigl (y,[\xi(y),u]_{\mathfrak m}\Bigr ) & = & 
- \sum_{i=1}^s \frac{C_i(y)}{C_k(y)} \cdot \alpha_i\Bigl (y,[y,u]_{\mathfrak m}\Bigr ) , \quad \forall u\in\fm.
\end{eqnarray}
Now, for any fixed $y\in\fm$, the values $C_i(y)$ are constants, $C_i(y)=c^i$ and determine
a naturally reductive metric from the initial family.
The components of the vector $\xi(y)$ are determined by the formula (\ref{ft2}) with $c^i=C_i(y)$
and consequently by the formula (\ref{ft3}).
Using the assumption that the formulas for the linear geodesic graph are equal for all the family
of the positively related naturally reductive Riemannian metric in question, the result follows.
$\hfill\square$

\begin{cor}
\label{c1}
Let $M=G/H$ be an irreducible homogeneous space with the ${\mathrm{Ad}}(H)$-irreducible decomposition
$\fm=\oplus_{i=1}^s\fm_i$ with at least two summands.
Let $g_{c^1\dots c^s}, c^i>0$, be a family of positively related naturally reductive Riemannian metrics on $M$
which admit unique geodesic graph.
Then there exist large family of Finsler metrics of the $\alpha_i$-type which are not naturally reductive.
\end{cor}
{\it{Proof.}}
Let $g_1,\dots,g_k$ be metrics from the family $g_{c^1\dots c^s}, c^i>0$.
We have observed that any Finsler metric $F=\sqrt{L(\sqrt{g_1},\dots,\sqrt{g_k})}$, for an admissible function $L$
(for example any metric of the type (\ref{Fex}) for $q\neq 2$ and $F_i=\sqrt{g_i}$),
is of the $\alpha_i$-type.
In general, the functions $k^j_r(C_1,\dots,C_s)$ do not cancel out to become real numbers
and we obtain unique nonlinear geodesic graph, which implies that the Finsler metric $F$ is g.o.,
but not naturally reductive.
$\hfill\square$

Let us remark that the assumption for the family of naturally reductive Riemannian metric to be positively
related is crucial for geodesic graph to be given by the formula (\ref{ft3}).
In \cite{AMD2}, we have described in detail a similar construction for initial g.o. metrics
on the sphere ${\mathrm{S}}^7$ which has decomposition $\fm=\fm_1\oplus\fm_2\oplus\fm_3$, it admits the 3-parameter
family of positively related g.o. metrics, but only the 2-parameter subfamily of these metrics is naturally reductive.
Metrics positively related with these
naturally reductive metrics are no more naturally reductive and for the construction
of the geodesic graph of the Finsler metric $F$ one has to use general formulas for nonlinear geodesic graph.
Hence, just for the family of positively related naturally reductive metrics as in Theorem above,
the geodesic graph given by formulas (\ref{ft3}) is the simplest possible, considering non Riemannian metrics
of the type $F=\sqrt{L(\sqrt{g_1},\dots,\sqrt{g_k})}$.

Let us consider now the same construction for the initial product manifold $M=M_1\times M_2$
(for example $M=H_3\times H_3$)
and the corresponding ${\mathrm{Ad}}(H)$-invariant decomposition $\fm=\fm_1\oplus\fm_2$
(where $\fm_i=\fn={\mathrm{Lie}}(H_3)$), which is not irreducible.
Here the basic scalar products $\alpha_i$ on $\fm_i$ determine invariant naturally reductive Riemannian metrics on $M_i$.
We consider again the family of metrics $g_{c^1,c^2}=\sum_{i=1}^2 c^i\alpha_i$.
In geodesic lemma given by the formula (\ref{f22}),
for given $u\in\fm_k$, on the left-hand side it holds $[\xi(y),u]\in\fm_k$ and also
on the right-hand side $[y,u]\in\fm_k$.
Hence the functions $C_i(y)$ cancel out of the formula.
Consequently, geodesic graph will be the sum of the linear geodesic graphs on $M_1$ and on $M_2$, respectively,
and the Finsler metric $F=\sqrt{L(\sqrt{g_1},\dots,\sqrt{g_k})}$ will be naturally reductive,
in accordance with the result of \cite{TX}.
This Finsler metric is the so-called $f$-product of the individual Riemannian naturally reductive metric on $M_i$,
for convenient function $f$ determined by $L$.
Our construction with any admissible function $L$ (for example the one given by formula (\ref{Fex}))
allows us to provide the explicit examples of Finsler naturally reductive metrics
which are $f$-products of two naturally reductive Riemannian metrics on the manifolds $M_i$.

\section{Finsler $F=\sqrt{L(\sqrt{g},\beta)}$ metrics}
\label{S4}
We shall now incorporate the one-forms into the Finsler metric.
First, we will consider Finsler metrics of the form $F=\sqrt{L(\sqrt{g},\beta)}$, which give alternative
description of the $(\alpha,\beta)$ metrics.
Assume that $g$ is a naturally reductive Riemannian metric.
The one-form $\beta$ must be also invariant.
We denote the corresponding ${\mathrm{Ad}}(H)$-invariant linear form on $\fm$ also by $\beta$
and we relate it with an ${\mathrm{Ad}}(H)$-invariant vector $v\in\fm$ using the formula
\begin{eqnarray}
\label{f24}
g(v,u) & = & \beta(u), \qquad \forall u\in\fm.
\end{eqnarray}
Consequently, the vector $v\in\fm$ determines an invariant Killing vector field $V$ on $M$.
In particular, the existence of an ${\mathrm{Ad}}(H)$-invariant vector $v\in\fm$ implies that at least one
of the subspaces in the ${\mathrm{Ad}}(H)$-irreducible decomposition $\fm=\oplus_{i=1}^s\fm_i$
must be one-dimensional, say $\fm_s={\mathrm{span}}(v)$.
First, let us consider locally irreducible manifolds, which do not admit parallel vector fields and hence
the one-form $\beta$ cannot be parallel.
\begin{prop}
\label{prop13}
Let $F=\sqrt{L(\sqrt{g},\beta)}$ be an invariant Finsler metric on $G/H$,
where $G$ is the maximal isometry group of $F$.
Let $(G/H,g)$ be irreducible Riemannian manifold which is
naturally reductive with respect to $G$ and with the unique geodesic graph with respect to $G$.
Then the geodesic graph of $F$ is nonlinear and the Finsler metric $F$ is not naturaly reductive.
\end{prop}
{\it{Proof.}}
The proof is based on the same idea as the proof of Propositions 4.1. and 4.4. in \cite{Du}.
We show first that there is the nontrivial center ${\mathfrak{c}}(\fg)\subset \fg$.
Consider the Killing form $B$ on $\fg$ and the $B$-orthogonal reductive decomposition $\fg=\fh+\fm_B$.
Denote by $v_B\in\fm_B$ the vector equivalent with $\beta\in T_p^*M$.
Assume $v_B\notin{\mathfrak{c}}(\fg)$.
Because the vector $v_B$ is ${\mathrm{Ad}}(H)$-invariant it holds $v_B\in C_\fg(\fh)$.
According to Lemma 8 in \cite{Ni}, for each $u\in C_\fg(\fh)$, it holds $[u,\fm_B]\subset \fm_B$
and ${\mathrm{ad}}(u)|_{\fm_B}$ is skew-symmetric operator.
In particular, ${\mathrm{ad}}(u)|_{\fm_B}$ is skew-symmetric for $u\in C_\fg(\fh)\cap\fm_B$.
Define the operator $w_B$ on $\fm_B$ (and consequently on $T_pM$) by $w_B(z)={\mathrm{ad}}(v_B)(z)$ for each $z\in\fm_B$.
Because $G$ is the maximal isometry group, we have $w_B\in\fh$ and $v_B-w_B\in C(\fg)$.
We now change the reductive decomposition for the naturally reductive decomposition $\fg=\fh+\fm$.
In particular, we consider the vector $v\in\fm$ which is the projection of the vector $v_B$ into $\fm$
and which again represents the one-form $\beta\in T_p^*M$.
If $v\in {\mathfrak{c}}(\fg)$,
the ${\mathrm{Ad}}(H)$-invariance of the vector $v$ and the natural reductivity leads to a contradiction
with the irreducibility of the manifold $M=G/H$.
Hence $v\notin {\mathfrak{c}}(\fg)$ and there exist a vector $w\in\fh$ such that $v-w\in {\mathfrak{c}}(\fg)$.
The geodesic lemma, which follows again from formula (\ref{fg}), now takes on the form
\begin{eqnarray}
\nonumber
\frac{ L_{,1}}{F_1(y)} \cdot g(y,[y+\xi(y),u]_\fm) + L_{,2}\cdot \beta([y+\xi(y),u])_\fm & = & 0, \quad \forall u\in\fm.
\end{eqnarray}
In the first term, we use the fact that the decomposition $\fg=\fh+\fm$ is naturally reductive.
In the second term, we use formula (\ref{f24}) and the ${\mathrm{Ad}}(H)$-invariance of the vector $v$.
We can rewrite the formula in the form
\begin{eqnarray}
\nonumber
g\Bigl (y,[\xi(y),u]_{\mathfrak m}\Bigr )
& = &
{F_1(y)}\cdot \frac{L_{,2}}{L_{,1}}
\cdot g \Bigl ( y, [ v ,u]_{\mathfrak m}\Bigr ), \quad \forall u\in\fm.
\end{eqnarray}
Now, we put
\begin{eqnarray}
\label{gg}
\xi(X)= {F_1(y)}\cdot \frac{L_{,2}}{L_{,1}} \cdot w
\end{eqnarray}
and the formula above is satisfied. Formula (\ref{gg}) describes the unique and nonlinear geodesic graph
and consequently the metric $F$ is not naturally reductive.
$~\hfill\square$

On the other hand, it is well known that the existence of a parallel vector field $V$ (or parallel one-form $\beta$)
implies the local reducibility of the manifold $M$.
In such a case, locally, the manifold $M=G/H$ can be expressed as $M=G_1/H\times{\mathbb{R}}$.
In the reductive decomposition $\fg=\fh+\fm$, one has $\fm=\fm_1\oplus\fm_2$,
where $\fm_1={\mathrm{Ker}}(\beta)$.
We do not have to relate the one-form $\beta$ with a vector $v$ and we get a result for a broader family
of positively related metrics.
\begin{prop}
\label{prop14}
Let $M_0=G/H$ be a homogeneous space with a reductive decomposition $\fg=\fm_0+\fh$.
Consider the product manifold $M=G/H\times{\mathbb{R}}$
with the ${\mathrm{Ad}}(H)$-irreducible decomposition $\fm=\oplus_{i=1}^{s} \fm_i=\fm_0\oplus\fm_s$,
positively related metrics $g_j$ with corresponding decompositions $g_j=\sum_{i=1}^s c_j^i \alpha_i$
and an invariant one-form $\beta$ such that $\fm_s={\mathrm{coKer}}(\beta)$.
The geodesic graph of the Finsler metric
$F_\beta=\sqrt{L(\sqrt{g_1},\dots,\sqrt{g_k},\beta)}$
coincides with the geodesic graph of the Finsler metric
$F_0=\sqrt{L(\sqrt{g_1},\dots,\sqrt{g_k},0)}$.
In particular, if $F_0$ is naturally reductive, then $F_\beta$ is also naturally reductive.
\end{prop}
{\it{Proof.}}
It holds $[\fm_0,\fm_s]=0$.
Because ${\mathrm{dim}}(\fm_s)=1$ and $\fm_s$ is an ${\mathrm{Ad}}(H)$-invariant subspace,
we obtain also $[\fg,\fm]_\fm\subset\fm_0={\mathrm{Ker}}(\beta)$.
Hence, the term involving $\beta$ in the geodesic lemma and formula (\ref{golemaf}) vanishes
and geodesic graph of the metric $F_\beta$ coincides with the geodesic graph for $F_0$.
$\hfill\square$

To conclude, let us relate these results with simple examples, in particular with
the Heisenberg group treated in Section \ref{sech}.
To illustrate Proposition \ref{prop13}, we consider any naturally reductive metric $g$ on $H_3$
and a one-form $\beta$ such that ${\mathrm{coKer}}(\beta)=\fm_2$ in the decomposition
$\fm=\fm_1\oplus\fm_2$ in Section \ref{sech}.
The Finsler metric $F=\sqrt{L(\sqrt{g},\beta)}$ can be considered as a metric of the $(\alpha_1,\alpha_2,\beta)$-type.
The geodesic graph for the metric $F$ can be described by formula (\ref{gg}) and the metric $F$ is not
naturally reductive.

To illustrate Proposition \ref{prop14},
we recall the classification result that any $4$-dimensional naturally reductive Riemannian manifold $M$
is locally a Riemannian product
$M=M_0\times {\mathbb{R}}$, where $M_0$ is a $3$-dimensional naturally reductive Riemannian manifold.
Consider, for example, the product manifold $M=H_3\times{\mathbb{R}}$.
On the Lie algebra level, we have the irreducible decomposition $\fm=\fm_1\oplus\fm_2\oplus\fm_3$.
Consider naturally reductive (and positively related) metrics $g_1,\dots,g_k$ on $M$
and an invariant one-form $\beta$ with ${\mathrm{coKer}}(\beta)=\fm_3$.
Metrics $F_\beta=\sqrt{L(\sqrt{g_1},\dots,\sqrt{g_k},\beta)}$
are of the $(\alpha_1,\alpha_2,\alpha_3,\beta)$-type and they have identical geodesic graph as the metrics
$F_0=\sqrt{L(\sqrt{g_1},\dots,\sqrt{g_k},0)}$.
If, in particular, each $g_j$ can be written as $k\cdot \alpha_0+l\cdot \alpha_3\, (k,l>0)$, where
$\alpha_0$ determines a naturally reductive metric on $H_3$ and
$\alpha_3$ determines a (naturally reductive) metric on ${\mathbb{R}}$,
then the metric $F_0$ is the $f$-product of $\alpha_0$ and $\alpha_3$
and both $F_0$ and $F_\beta$ are naturally reductive.
In general, our construction allows to consider two independent and invariant one-forms $\beta_1, \beta_2$ on $M$,
namely those determined by the conditions ${\mathrm{coKer}}(\beta_1)=\fm_2$ and ${\mathrm{coKer}}(\beta_2)=\fm_3$,
and Finsler metrics $F_{\beta_1,\beta_2}=\sqrt{L(\sqrt{g_1},\dots,\sqrt{g_k},\beta_1,\beta_2)}$
lead to explicit examples of metrics of the type $(\alpha_1,\alpha_2,\alpha_3,\beta_1,\beta_2)$.
Metrics $F_{\beta_1,\beta_2}=\sqrt{L(\sqrt{g_1},\dots,\sqrt{g_k},\beta_1,\beta_2)}$
and $F_{\beta_1,0}=\sqrt{L(\sqrt{g_1},\dots,\sqrt{g_k},\beta_1,0)}$ have identical geodesic graphs.

\section*{Statements and Declarations}

\textbf{Funding:}
{The research is supported by grant GR24068
funded by Junta de Extremadura, partially funded by Fondo Europeo de Desarrollo Regional.}

\noindent
\textbf{Competing Interests:}
The authors have no financial or competing interests to declare that are relevant to the content of this article.

\noindent
\textbf{Author Contributions:}
Both authors contributed equally to this research and in writing the paper.

\noindent
\textbf{Availability of data and material:}
Not applicable.


\begin{thebibliography}{[0000]}
\bibitem {AMD}
Arias-Marco, T., Du\v sek, Z.:
Geodesic graphs for geodesic orbit Finsler $(\alpha,\beta)$ metrics on spheres,
{\it Vietnam J. Math.} (2024). https://doi.org/10.1007/s10013-024-00695-x
\bibitem {AMD2}
Arias-Marco, T., Du\v sek, Z.:
Structure of geodesics for Finsler metrics arising from the Riemannian g.o. metrics,
{\it Kyungpook J. Math.}, to appear.
\bibitem {BN}
Berestovskii, V., Nikonorov, Yu.:
Riemannian Manifolds and Homogeneous Geodesics,
Springer Nature Switzerland, Cham, 2020.
\bibitem {DH}
Deng, S., Hou, Z.:
Invariant Finsler metrics on homogeneous manifolds,
{\it J. Phys. A: Math. Gen.} {\bf{37}} (2004), 8245--8253.
\bibitem {DX}
Deng, S., Xu, M.:
$(\alpha_1,\alpha_2)$-Metrics and Clifforf-Wolf Homogeneity,
{\it J. Geom. Anal.} {\bf{26}} (2016), 2282--2321.
\bibitem {DYZ}
Deng, S., Yan, Z., Zhang, S.:
Naturally reductive homogeneous $(\alpha_1,\alpha_2)$ spaces,
{\it Publ. Math. Debrecen} {\bf{102/3-4}} (2023), 415--427.
\bibitem{DuS2}
Du\v sek, Z.:
Homogeneous geodesics and g.o. manifolds,
{\it Note Mat.} {\bf{38}} (2018), 1--15.
\bibitem{DuCMUC}
Du\v sek, Z.:
Geodesic graphs in Randers g.o. spaces,
{\it Comment. Math. Univ. Carol.} {\bf{61,2}} (2020), 195--211.
\bibitem{DuWS}
Du\v sek, Z.:
Structure of geodesics in weakly symmetric Finsler metrics on H-type groups,
{\it Arch. Math.-Brno} {\bf{56}} (2020), 265--275.
\bibitem{Du}
Du\v sek, Z.:
Geodesic orbit Finsler $(\alpha,\beta)$ metrics,
{\it Eur. J. Math.} {\bf{9,1}} (2023), 9.
\bibitem {JS}
Javaloyes, M.A., S\'anchez, M.:
On the definition and examples of Finsler metrics,
{\it Ann. Scuola Norm. Super. Pisa-Cl. Sci.} Vol. XIII (2014), 813--858.
\bibitem {KNi}
Kowalski, O., Nik\v cevi\' c, S.:
On geodesic graphs of Riemannian g.o. spaces,
{\it Arch. Math.} {\bf{73}} (1999), 223--234;
Appendix:
{\it Arch. Math.} \textbf{79} (2002), 158--160.
\bibitem {La}
Latifi, D.:
Homogeneous geodesics in homogeneous Finsler spaces,
{\it J. Geom. Phys.} {\bf{57}} (2007), 1421--1433.
\bibitem {Ni}
Nikonorov, Yu.G.:
On the structure of geodesic orbit Riemannian spaces,
{\it Ann. Glob. Anal. Geom.} {\bf{52}}(3) (2017), 289--311.
\bibitem {Sz} Szenthe, J.:
Sur la connection naturelle \`a torsion nulle,
{\it Acta Sci. Math. (Szeged)} {\bf{38}} (1976), 383--398.
\bibitem {TX}
Tan, J., Xu, M.:
Naturally reductive $(\alpha_1,\alpha_2)$ metrics
{\it Acta Math. Sci.} {\bf{43B}}(4) (2023), 1547--1560.
\bibitem {Xu}
Xu, M.:
Geodesic orbit spheres and constant curvature in Finsler geometry,
{\it Differ. Geom. Appl.} {\bf{61}} (2018), 197--206.
\bibitem {Y}
Yan, Z.:
Some Finsler spaces with homogeneous geodesics,
{\it Math. Nachr.} {\bf{290}},2-3 (2017), 474--481.
\bibitem {YD}
Yan, Z., Deng, S.:
Finsler spaces whose geodesics are orbits,
{\it Differ. Geom. Appl.} {\bf{36}} (2014), 1--23.
\end{thebibliography}
\end{document}